\documentclass[12pt]{article}
\usepackage{amsmath}
\usepackage{amssymb}
\usepackage{amscd}
\usepackage{amsthm}

\theoremstyle{plain}
\newtheorem{Thm}{Theorem}[section]

\newtheorem{Cor}[Thm]{Corollary}
\newtheorem{Lem}[Thm]{Lemma}

\theoremstyle{definition}

\newtheorem{Expl}[Thm]{Example}

\newtheorem{Rem}[Thm]{Remark}

\numberwithin{equation}{section}

\title{Derived Categories of Toric Varieties}
\author{Yujiro Kawamata}

\begin{document}

\maketitle

\section{Introduction}

The purpose of this paper is to investigate the structure of the derived 
category of a toric variety.
We shall prove the following:

\begin{Thm}\label{main}
Let $X$ be a projective toric variety with at most quotient singularities, 
let $B$ be an invariant $\mathbb{Q}$-divisor whose coefficients belong to the
set $\{\frac{r-1}r ; r \in \mathbb{Z}_{>0}\}$, and let $\mathcal{X}$ be the 
smooth Deligne-Mumford stack associated to the pair $(X,B)$ as in 
\cite{log-crep}.
Then the bounded derived category of coherent sheaves
$D^b(\text{Coh}(\mathcal{X}))$ has a complete exceptional collection 
consisting of sheaves.
\end{Thm}

An object of a triangulated category $a \in T$ is called {\em exceptional}
if 
\[
\text{Hom}^p(e,e) \cong \begin{cases} \mathbb{C} &\text{ for } p = 0 \\
0 &\text{ for } p \ne 0. \end{cases}
\]
A sequence of exceptional objects $\{e_1,\dots,e_m\}$ is said to be an 
{\em exceptional collection} if 
\[
\text{Hom}^p(e_i,e_j) = 0 \text{ for all } p \text{ and } i > j.
\]
It is said to be {\em strong} if in addition that 
$\text{Hom}^p(e_i,e_j) = 0$ for $p \ne 0$ and all $i,j$. 
It is called {\em complete} if $T$ coincides with the smallest triangulated
subcategory containing all the $e_i$ (cf. \cite{Bondal}).

It is usually hard to determine the explicit structure of a derived 
category of a variety.
But it is known that some special varieties such as a projective space or
a Grassmann variety have strong complete exceptional collections consisting of 
vector bundles
(\cite{Beilinson}, \cite{Kapranov1}, \cite{Kapranov2}, \cite{Kapranov3}).
Such sheaves are useful for further investigation of the derived categories
(\cite{Rudakov}, \cite{GKR}, \cite{Macri} and \cite{Bridgeland2} for example).

We use the minimal model program for toric varieties as developed in 
\cite{Reid} (and corrected in \cite{Matsuki}) in order to prove the theorem.
A special feature for this approach is that, 
even if we start with the smooth and non-boundary case $B = 0$, 
we are forced to deal with not only singularities but also the case $B \ne 0$ 
because Mori fiber spaces have multiple fibers in general.
Thus we are inevitably lead to consider
the general situation concerning Deligne-Mumford stacks
even if we only need results for smooth varieties.
The stacky sheaves need careful treatment 
because there exist non-trivial stabilizer groups on the stacks
(cf. Remark~\ref{stabilizer}), 

We start with the Beilinson theorem for the case of projective spaces, and
build up exceptional collections following the procedure of 
the minimal model program.
We use a covering trick to proceed from projective spaces to log Fano
varieties (\S 3).
Then we proceed by induction on the dimension.
First we consider a Mori fiber space in \S 4 
where the base space is assumed to have
already a complete exceptional collection by the induction hypothesis.
Though a Mori fiber space has singular fibers, the associated morphism of 
stacks is proved to be smooth (Corollary~\ref{smooth}), and we can 
define a complete exceptional collection on the total space by using twisted 
pull-backs.
The behavior of derived catogories under the birational transformations such as
a divisorial contraction or a flip was studied in \cite{log-crep}.
We use this result together with results of \S 4 in \S 5 and \S 6.
Indeed, the exceptional locus of a divisorial contraction or a flip
has a structure of a Mori fiber space itself.
The argument of the proof is a generalization of that in \cite{Orlov} 
which considered the derived categories of a projective space bundle 
and a blowings-up with a smooth center in a smooth variety.


\section{Toric minimal model program}

Let $X$ be a projective toric variety of dimension $n$ which is quasi-smooth,
i.e., has only quotient singularities.
We note that a toric variety is quasi-smooth if and only if it is 
$\mathbb{Q}$-factorial.
We consider a $\mathbb{Q}$-divisor $B$ on $X$ whose prime components are 
invariant divisors with coefficients being contained in the set 
$\{\frac{r-1}r ; r \in \mathbb{Z}_{>0}\}$.
Let $\mathcal{X}$ be the smooth Deligne-Mumford stack associated 
to the pair $(X,B)$ 
with the natural morphism $\pi_X: \mathcal{X} \to X$ as in 
\cite{log-crep}.

The pair $(X,B)$ has only log terminal singularities.
We work on the log minimal model program for $(X,B)$.
We refer the reader to \cite{Reid} or \cite{Matsuki}. 
Let $\phi: X \to Y$ be a primitive contraction morphism corresponding to an
extremal ray with respect to $K_X+B$.
Then $Y$ is also a projective toric variety and $\phi$ is a toric morphism.
If $\phi$ is a birational morphism, then the boundary divisor $C$ on $Y$ 
is defined to be the strict transform of $B$.
Otherwise, it will be defined later.

Let $N_X$ be the lattice of $1$-parameter subgroups of the torus acting on 
$X$, and $\Delta_X$ the fan in $N_{X,\mathbb{R}}$ corresponding to $X$.
Let $w = \langle v_3, \dots, v_{n+1} \rangle$ 
be a wall in $\Delta_X$ corresponding to an extremal rational curve, where 
the $v_i$ are primitive vectors in $N_X$ on the edges of $w$.
Let $v_1$ and $v_2$ be two primitive vectors in $N_X$ each of which forms 
an $n$-dimensional cone in $\Delta_X$ when combined with $w$. 
Let $D_i$ be the prime divisors on $X$ corresponding to the $v_i$, and 
$\mathcal{D}_i$ the corresponding prime divisors on $\mathcal{X}$.
Let $\frac{r_i-1}{r_i}$ be the coefficients of the $D_i$ in $B$.
Then the natural morphism $\pi_X: \mathcal{X} \to X$ ramifies along
$D_i$ such that $\pi_X^*D_i = r_i\mathcal{D}_i$.

The contraction morphism is described by an equation
\begin{equation}\label{relation}
a_1v_1 + \dots + a_{n+1}v_{n+1} = 0
\end{equation}
where the $a_i$ are integers such that 
\[
\begin{split}
&(a_1, \dots, a_{n+1})=1 \\
&a_i > 0 \text{ for } 1 \le i \le \alpha \\
&a_i = 0 \text{ for } \alpha + 1 \le i \le \beta \\
&a_i < 0 \text{ for } \beta + 1 \le i \le n + 1 \\
&2 \le \alpha \le \beta \le n+1.
\end{split}
\]
Note that we use slightly different notation from the literatures
where the $a_i$ are rational numbers.
Since $K_X + B$ is negative for $\phi$, 
we have 
\[
\sum_{i=1}^{n+1} \frac{a_i}{r_i} > 0.
\]
The following lemma asserts that the 
set of integers $\{a_i\}$ is {\em well prepared}:

\begin{Lem}\label{prepared}
Let $i_0$ be an integer such that $1 \le i_0 \le \alpha$ or 
$\beta + 1 \le i_0 \le n + 1$.
Then the set of $n+\alpha-\beta$ integers $a_i$ for $1 \le i \le \alpha$ and 
$\beta + 1 \le i \le n + 1$ except one $i = i_0$ is coprime for any $i_0$.
\end{Lem}

\begin{proof}
Let $c$ be the largest common divisor of these integers, and 
set $a_i = c\bar a_i$.
Then we have $(a_{i_0},c) = 1$.
Let $x,y$ be integers such that $a_{i_0}x + cy = 1$.
Then 
\[
\frac 1c v_{i_0} = \frac 1c v_{i_0} - \frac xc \sum_i a_iv_i
= - x \sum_{i \ne i_0} \bar a_i v_i + yv_{i_0} \in N_X
\]
hence $c = 1$.
\end{proof}

One of the following cases occurs.

(1) Mori fiber space: $\beta = n+1$.
We have $\dim Y = n+1-\alpha$.

(2) divisorial contraction: $\beta = n$.

(3) small contraction: $\beta < n$. 

We treat these cases separately in the following sections.
According to the minimal model program, Theorem~\ref{main} follows from the 
combination of Corollaries~\ref{cfiber}, \ref{cdivisorial} and \ref{cflip}.


\section{Fano case}

We start with the case where $(X,B)$ is a log $\mathbb{Q}$-Fano variety 
with $\rho = 1$.
We have $\alpha = \beta = n+1$.
In this case, there are no edges in $\Delta_X$ besides $\mathbb{R}_{\ge 0}v_i$.
Such a variety $X$ is not necessarily a weighted projective space
as remarked in \cite{Matsuki}.
But it is covered by a weighted projective space by a finite morphism which
is etale in codimension $1$.
Indeed, a weighted projective space is characterized by the property that 
the divisor class group has no torsion (cf. Lemma~\ref{torsion}).

Let $N_{X'}$ be the sublattice of $N_X$ generated by the $v_i$.
By the equation (\ref{relation}), 
the toric variety $X'$ corresponding to the fan $\Delta_X$ in
$N_{X',\mathbb{R}}$, with the lattice $N_{X'}$, 
is isomorphic to the weighted projective space
$\mathbb{P}(a_1,\dots,a_{n+1})$.
The natural morphism $\bar{\sigma}_1: X' \to X$ is etale in codimension $1$.
Let $D'_i$ be the prime divisors on $X'$ corresponding to the $v_i$, 
$\mathcal{X}'$ the smooth Deligne-Mumford stack associated to the pair 
$(X', \sum_i \frac{r_i-1}{r_i}D'_i)$ with the projections
$\pi_{X'}: \mathcal{X}' \to X'$ and $\sigma_1: \mathcal{X}' \to \mathcal{X}$, 
and let $\mathcal{D}'_i$ be the prime divisors on 
$\mathcal{X}'$ such that $\pi_{X'}^*D'_i = r_i\mathcal{D}'_i$.

Let $r$ be a positive integer such that $a_ir$ is divisible by $r_i$ for 
any $i$, and let $N_{\tilde X}$ be the sublattice of $N_X$ generated by the
vectors $\tilde v_i = a_irv_i$.
We have
\[
\sum_{i=1}^{n+1} \tilde v_i = 0,
\]
and the toric variety $\tilde X$ corresponding to the fan $\Delta_X$ in
$N_{\tilde X,\mathbb{R}}$, with the lattice $N_{\tilde X}$, 
is isomorphic to the projective space $\mathbb{P}^n$.
Let $\sigma_2: \tilde X \to \mathcal{X}'$ be the natural morphism, and set 
$\sigma = \sigma_1 \circ \sigma_2$, $\bar{\sigma}_2 = \pi_{X'} \circ \sigma_2$
and $\bar{\sigma} = \pi_X \circ \sigma$. 
Let $\tilde D_i$ be the prime divisors on $\tilde X$ corresponding to the 
vectors $\tilde v_i$.
Moreover, let $N_{X''}$ be the sublattice of $N_X$ generated by the $r_iv_i$. 
We note that 
the vectors $r_iv_i$ are not necessarily primitive in this lattice.

\begin{Lem}\label{torsion}
(1) A divisor $\sum_i k_i \mathcal{D}_i$ is torsion 
in the divisor class group of $\mathcal{X}$ if 
and only if 
\[
\sum_i \frac{a_ik_i}{r_i} = 0.
\]

(2) The group of torsion divisor 
classes on $\mathcal{X}$ is dual to the quotient group $N_X/N_{X''}$. 

(3) The group of torsion Weil divisor 
classes on $X$ is dual to the quotient group $N_X/N_{X'}$. 
\end{Lem}

\begin{proof}
(1) A divisor $\sum_i k_i \mathcal{D}_i$ is linearly equivalent to $0$ 
if and only if 
there exists $m \in M_X = N_X^*$ such that $(m,r_iv_i) = k_i$,
because the morphism $\pi_X: \mathcal{X} \to X$ is birational.
Thus $\sum_i k_i \mathcal{D}_i$ is torsion if and only if 
there exists $m \in M_{X,\mathbb{R}}$ such that $(m,r_iv_i) = k_i$.
The latter condition is equivalent to the equality 
$\sum_i \frac{a_ik_i}{r_i} = 0$.

(2) For $m \in M_{X,\mathbb{R}}$, we have $(m,r_iv_i) \in \mathbb{Z}$ 
for all $v_i$
if and only if $m \in M_{X''} = N_{X''}^*$.
Therefore, the group of torsion divisor classes is isomorphic to 
$M_{X''}/M_X$.

(3) is a particular case of (2).
\end{proof}

\begin{Rem}
If $B=0$, i.e., $r_i=1$ for all $i$, then 
the divisor class groups of $X$ and $\mathcal{X}$ are 
isomorphic.
\end{Rem}

\begin{Expl}
Torsion divisor classes correspond to etale coverings of the stack.
For example, let $X = \mathbb{P}^n$ be the projective space, and
$\mathcal{X}$ the smooth stack associated to the pair
$(X, \frac{r-1}r \sum_{i=1}^{n+1} H_i)$, where the $H_i$ are 
coordinate hyperplanes.
Let $\mathcal{H}_i$ be the prime divisors on $\mathcal{X}$ above the $H_i$
so that $\pi_X^*H_i = r\mathcal{H}_i$ for the projection 
$\pi_X: \mathcal{X} \to X$.

Let $\tilde X = \mathbb{P}^n$ be another projective space, 
and let $\sigma: \tilde X \to \mathcal{X}$
be the Kummer covering with Galois group $(\mathbb{Z}/r)^n$ obtained by 
taking the $r$-th roots of the coordinates.
Then $\sigma$ is etale, and we have
\[
\sigma_*\mathcal{O}_{\tilde X}
\cong \bigoplus_{l_1,\dots,l_n = 0}^{r-1} \mathcal{O}_{\mathcal{X}}
(\sum_{i=1}^n l_i \mathcal{H}_i +(- \sum_{i=1}^n l_i)\mathcal{H}_{n+1}).
\]
We note that the direct summands are invertible sheaves on $\mathcal{X}$
corresponding to the torsion divisor classes.
In the usual language, if we denote 
$\bar{\sigma} = \pi_X \circ \sigma: \tilde X \to X$,
then
\[
\bar{\sigma}_*\mathcal{O}_{\tilde X}
\cong \bigoplus_{l_1,\dots,l_n = 0}^{r-1} \mathcal{O}_X
(- \ulcorner \frac{\sum_{i=1}^n l_i}r \urcorner H_{n+1})
\]
because
\[
\pi_{X*}\mathcal{O}_{\mathcal{X}}(\sum_{i=1}^n l_i \mathcal{H}_i 
+(- \sum_{i=1}^n l_i)\mathcal{H}_{n+1})
= \mathcal{O}_X(- \ulcorner \frac{\sum_{i=1}^n l_i}r \urcorner H_{n+1}).
\]
More generally, we have 
\[
\sigma_*\mathcal{O}_{\tilde X}(-p)
\cong \bigoplus_{l_1,\dots,l_n = 0}^{r-1} \mathcal{O}_{\mathcal{X}}
(\sum_{i=1}^n l_i \mathcal{H}_i 
+ (- p - \sum_{i=1}^n l_i)\mathcal{H}_{n+1})
\]
and
\[
\bar{\sigma}_*\mathcal{O}_{\tilde X}(-p)
\cong \bigoplus_{l_1,\dots,l_n = 0}^{r-1} \mathcal{O}_X
(- \ulcorner \frac{p+\sum_{i=1}^n l_i}r \urcorner H_{n+1}).
\]
For example, the direct images of the sheaves $\mathcal{O}_{\tilde X}(-p)$
for $0 \le p \le n$, which generates the derived category 
$D^b(\text{Coh}(\tilde X))$ (\cite{Beilinson}), 
have the direct summands of the form 
$\mathcal{O}_X(- q)$ for $0 \le q \le n$
(cf. \cite{equi}).
\end{Expl}

\begin{Lem}\label{decomp}
(1) Let $G_1 = N_X/N_{X'}$ be the Galois group of the covering 
$\bar{\sigma}_1: X' \to X$.
Then there is the following decomposition into eigenspaces with respect to the 
$G_1$-action
\[
\sigma_{1*}\mathcal{O}_{\mathcal{X}'}(\sum_i d_i \mathcal{D}'_i)
\cong \bigoplus_k \mathcal{O}_{\mathcal{X}}(\sum_i (d_i+k_ir_i)\mathcal{D}_i)
\]
where the sequences of integers $k = (k_i)$ 
in the summation are determined by the equation $k_i = (m,v_i)$ for 
the representatives $m$ 
of the group of torsion Weil divisor classes $M_{X'}/M_X$ of $X$.

(2) Let $G_2 = N_{X'}/N_{\tilde X}$ be the Galois group of the covering 
$\bar{\sigma}_2: \tilde X \to X'$.
Then there is the following decomposition into eigenspaces with respect to the 
$G_2$-action
\[
\sigma_{1*}\mathcal{O}_{\tilde X}(-p) 
\cong \bigoplus_l \mathcal{O}_{\mathcal{X}'}
(\sum_{1 \le i \le n} \llcorner \frac{l_ir_i}{a_ir} \lrcorner \mathcal{D}'_i
+ \llcorner \frac{(l_{n+1}-p)r_{n+1}}{a_{n+1}r} \lrcorner \mathcal{D}'_{n+1}
- \frac 1r \sum_{i=1}^{n+1} l_i)
\]
where the sequences of integers $l = (l_i)$ in the summation
run under the conditions that $0 \le l_i < a_ir$ and 
$r \vert \sum_{i=1}^{n+1} l_i$.

(3) Let $G = N_X/N_{\tilde X}$ the Galois group of the covering 
$\bar{\sigma}: \tilde X \to X$.
Then there is the following decomposition into eigenspaces with respect to the 
$G$-action
\[
\sigma_*\mathcal{O}_{\tilde X}(-p) 
\cong \bigoplus_k \mathcal{O}_{\mathcal{X}}
(\sum_{1 \le i \le n} \llcorner \frac{k_ir_i}{a_ir} \lrcorner \mathcal{D}_i
+ \llcorner \frac{(k_{n+1}-p)r_{n+1}}{a_{n+1}r} \lrcorner \mathcal{D}_{n+1})
\]
where the sequences of integers $k = (k_i)$ satisfy the equation 
\[
\sum_{i=1}^{n+1} k_i = 0.
\]
\end{Lem}

\begin{proof}
(1) is clear.

(2) We have an exact sequence
\[
0 \to \mathbb{Z}/r \to \bigoplus_{i=1}^{n+1} \mathbb{Z}/a_ir \to G_2 \to 0
\]
where $1$ in the first term is sent to $(a_i)$ in the second term.
Thus 
\[
G_2^* \cong \{(l_i) \in \bigoplus_{i=1}^{n+1} \mathbb{Z}/a_ir ; 
\sum_{i=1}^{n+1} l_i = 0 \mod r \}.
\]
We have $\sigma_1^*\mathcal{D}'_i = \frac{a_ir}{r_i} \tilde D_i$
for the prime divisor $\tilde D_i$ on $\tilde X$ above $\mathcal{D}'_i$.
Since $\frac{a_il_i}{a_ir}=\frac{l_i}r$, we obtain the formula. 
We note that $\mathcal{O}_{\mathcal{X}'}(1)$ is well-defined because 
$\mathcal{X}'$ has no torsion divisor classes.

(3) By combining (1) and (2), we obtain
\[
\begin{split}
\sigma_*\mathcal{O}_{\tilde X}(-p) 
\cong \bigoplus_{k,l} \mathcal{O}_{\mathcal{X}}
(&\sum_{1 \le i \le n} (k_ir_i + \llcorner \frac{l_ir_i}{a_ir} \lrcorner)
\mathcal{D}_i \\
&+ (k_{n+1}r_{n+1} + \llcorner \frac{(l_{n+1}-p)r_{n+1}}{a_{n+1}r} \lrcorner) 
\mathcal{D}_{n+1})
\end{split}
\]
where the $k = (k_i)$ satisfy the equation 
\[
\sum_{i=1}^{n+1} a_ik_i = - \frac 1r \sum_{i=1}^{n+1} l_i
\]
and the summation on $l = (l_i)$ 
is under the restriction that $0 \le l_i < a_ir$ and $r \vert \sum_i l_i$.
If we replace $a_irk_i + l_i$ by $k_i$, then we obtain our assertion.
\end{proof}

\begin{Thm}\label{Fano}
(1) An invertible sheaf 
$\mathcal{O}_{\mathcal{X}}(\sum_{i=1}^{n+1} k_i \mathcal{D}_i)$ 
on $\mathcal{X}$ is an exceptional object for any sequence of integers 
$k = (k_i)$ for $1 \le i \le n+1$.

(2) If $\sum_{i=1}^{n+1} \frac{a_ik_i}{r_i} 
> \sum_{i=1}^{n+1} \frac{a_ik'_i}{r_i} > 
\sum_{i=1}^{n+1} \frac{a_i(k_i-1)}{r_i}$, then
\[
\text{Hom}^q(\mathcal{O}_{\mathcal{X}}(\sum_{i=1}^{n+1} k_i \mathcal{D}_i),
\mathcal{O}_{\mathcal{X}}(\sum_{i=1}^{n+1} k'_i \mathcal{D}_i)) = 0
\]
for all $q$, where $k' = (k'_i)$ is another sequence of integers.

(3) If $\sum_{i=1}^{n+1} \frac{a_ik_i}{r_i} 
= \sum_{i=1}^{n+1} \frac{a_ik'_i}{r_i}$ and
$\sum_{i=1}^{n+1} k_i \mathcal{D}_i \not\sim 
\sum_{i=1}^{n+1} k'_i \mathcal{D}_i$, then
\[
\text{Hom}^q(\mathcal{O}_{\mathcal{X}}(\sum_{i=1}^{n+1} k_i \mathcal{D}_i),
\mathcal{O}_{\mathcal{X}}(\sum_{i=1}^{n+1} k'_i \mathcal{D}_i)) = 0
\]
for all $q$.

(4) If $\sum_{i=1}^{n+1} \frac{a_ik_i}{r_i} 
\le \sum_{i=1}^{n+1} \frac{a_ik'_i}{r_i}$, then
\[
\text{Hom}^q(\mathcal{O}_{\mathcal{X}}(\sum_{i=1}^{n+1} k_i \mathcal{D}_i),
\mathcal{O}_{\mathcal{X}}(\sum_{i=1}^{n+1} k'_i \mathcal{D}_i)) = 0
\]
for $q \ne 0$.

(5) The set of invertible sheaves 
$\mathcal{O}_{\mathcal{X}}(\sum_{i=1}^{n+1} k_i \mathcal{D}_i)$ for 
\[
0 \ge \sum_{i=1}^{n+1} \frac{a_ik_i}{r_i} > - \sum_{i=1}^{n+1} \frac{a_i}{r_i}
\]
generates the triangulated category $D^b(\text{Coh}(\mathcal{X}))$.
\end{Thm}

\begin{proof}
The canonical divisor of $\mathcal{X}$ is given by 
\[
\omega_{\mathcal{X}} \cong \pi_X^*\omega_X \otimes 
\mathcal{O}_{\mathcal{X}}(\sum_i (r_i-1)\mathcal{D}_i)
\cong \mathcal{O}_{\mathcal{X}}(- \sum_i \mathcal{D}_i).
\] 
An invertible sheaf $\mathcal{O}_{\mathcal{X}}(\sum_i k_i \mathcal{D}_i)$
is ample if and only if $\sum_i \frac{a_ik_i}{r_i} > 0$.
Therefore, the assertions (1) through (4) follow immediately 
from the vanishing theorem (\cite{KMM}).

(5) follows from a similar generalization of the Beilinson resolution 
theorem (\cite{Beilinson}) as in \cite{equi}~\S 5.
Indeed, the integral functor corresponding to an object $e$ on 
$\mathcal{X} \times \mathcal{X}$ given by
\[
\begin{split}
e = &\{0 \to 
[\sigma_*\mathcal{O}_{\tilde X}(-n) \boxtimes 
\sigma_*\Omega^n_{\tilde X}(n)]^G \to \dots \\
&\to [\sigma_*\mathcal{O}_{\tilde X}(-1) \boxtimes 
\sigma_*\Omega^1_{\tilde X}(1)]^G \to 
[\sigma_*\mathcal{O}_{\tilde X} \boxtimes 
\sigma_*\mathcal{O}_{\tilde X}]^G \to 0\}
\end{split}
\]
is isomorphic to the identity functor, where the group $G$ acts diagonally
on the tensor products.
Thus the derived category $D^b(\text{Coh}(\mathcal{X}))$ is generated by the
direct summands of the sheaves $\sigma_*\mathcal{O}_{\tilde X}(-p)$ for 
$0 \le p \le n$ given in Lemma~\ref{decomp}~(3).

Since $\sum_i k_i = 0$, we have 
\[
\sum_{i=1}^n \frac{a_i\frac{k_ir_i}{a_ir}}{r_i}+
\frac{a_{n+1}\frac{(k_{n+1}-p)r_{n+1}}{a_{n+1}r}}{r_{n+1}}
= - \frac pr.
\]
Then we calculate
\[
\begin{split}
&0 \ge \sum_{i=1}^n \frac{a_i\llcorner \frac{k_ir_i}{a_ir} \lrcorner}{r_i}+
\frac{a_{n+1}\llcorner \frac{(k_{n+1}-p)r_{n+1}}{a_{n+1}r} \lrcorner}{r_{n+1}} 
\\
&\ge - \sum_{i=1}^{n+1} \frac{a_i}{r_i}(1 - \frac 1{s_i}) - \frac nr
\ge - \sum_{i=1}^{n+1} \frac{a_i}{r_i} + \frac 1r
\end{split}
\]
where we put $a_ir = r_is_i$ for some integers $s_i$.
\end{proof}

\begin{Cor}\label{Fano}
Let $(X,B)$ be a $\mathbb{Q}$-factorial projective 
toric variety such that $- (K_X+B)$ is ample, $\rho(X) = 1$ and that 
the coefficients of $B$ belong 
to the set $\{\frac{r-1}r ; r \in \mathbb{Z}_{>0}\}$.
Let $\mathcal{X}$ be the smooth Deligne-Mumford stack associated to the 
pair $(X,B)$.
Then the derived category 
$D^b(\text{Coh}(\mathcal{X}))$ has a strong complete exceptional 
collection consisting of invertible sheaves.
\end{Cor}

\begin{proof}
The number of isomorphism classes of the set of invertible sheaves 
$\mathcal{O}_{\mathcal{X}}(\sum_i k_i \mathcal{D}_i)$ for 
$0 \ge \sum_i \frac{a_ik_i}{r_i} > - \sum_i \frac{a_i}{r_i}$ is finite.
\end{proof}


\section{Mori fiber space}

We consider a toric Mori fiber space $\phi: X \to Y$ with respect to 
$K_X+B$.
This fibration is not necessarily locally trivial because there may be 
multiple fibers.
But it becomes locally trivial after taking coverings:

\begin{Lem}
Let $Y_0$ be an invariant open affine subset of $Y$, 
and let $X_0 = \phi^{-1}(Y_0)$.
Then there exist finite surjective toric morphisms 
$\tau_{X_0}: X'_0 \to X_0$ and 
$\tau_{Y_0}: Y''_0 \to Y_0$ with a toric surjective morphism 
$\phi'_0: X'_0 \to Y''_0$ 
which satisfy the following conditions.

(1) $\tau_{X_0}$ is etale in codimension $1$.

(2) $\phi \circ \tau_{X_0} = \tau_{Y_0} \circ \phi'_0$.

(3) $X'_0$ is isomorphic to the direct product of $Y''_0$ and a weighted 
projective space, and
$\phi'_0$ corresponds to the projection.
\end{Lem}

\begin{proof}
Let $N_Y$ be the lattice of $1$-parameter subgroups of the torus for $Y$, and 
$\Delta_Y$ the fan in $N_{Y,\mathbb{R}}$ corresponding to $Y$.
We take the wall $w$ described in the formula (\ref{relation}) such that 
the corresponding extremal rational curve is contained in $X_0$.
We have 
\[
N_Y = N_X/(\bigoplus_{i=1}^{\alpha} \mathbb{R}v_i \cap N_X).
\]
Let $h: N_X \to N_Y$ be the projection.
We write $h(v_i) = s_i\bar v_i$ for primitive vectors $\bar v_i$ in $N_Y$ 
and positive integers $s_i$ for $\alpha+1 \le i \le n+1$.
Then these $\bar v_i$ give the set of edges of an $(n+1-\alpha)$-dimensional 
cone $\sigma_0$ in $\Delta_Y$ corresponding to $Y_0$.
Let $E_i$ be the prime divisors on $Y$ corresponding to the vectors $\bar v_i$.
$X_0$ coincides with the toric variety corresponding to the fan 
$\Delta_X \cap h^{-1}(\sigma_0)$ in $N_{X,\mathbb{R}}$.

Let $N_{X'_0}$ be the sublattice of $N_X$ 
generated by the $v_i$ for $1 \le i \le n+1$, and
$N_{Y'_0}$ (resp. $N_{Y''_0}$) of $N_Y$
generated by the $\bar v_i$ (resp. $h(v_i)$)
for $\alpha+1 \le i \le n+1$. 
Let $X'_0$ be the toric variety corresponding to the fan 
$\Delta_X \cap h^{-1}(\sigma_0)$ in $N_{X'_0,\mathbb{R}}$, and 
$Y'_0$ (resp. $Y''_0$) the one corresponding to the 
cone $\sigma_0$ in $N_{Y'_0,\mathbb{R}}$ 
(resp. $N_{Y''_0,\mathbb{R}}$).
Then the natural morphisms $\tau_{X_0}: X'_0 \to X$ and 
$\tau'_{Y_0}: Y'_0 \to Y$ are etale in codimension $1$, 
while $\tau''_{Y_0}: Y''_0 \to Y'_0$ is not in general.
Since 
\[
\sum_{i=1}^{\alpha} a_iv_i = 0
\]
$X'_0$ is isomorphic to the product of $Y''_0$ 
with a weighted projective space $\mathbb{P}(a_1,\dots,a_{\alpha})$.
\end{proof}

We define the boundary $\mathbb{Q}$-divisor $C$ on $Y$ by assigning 
coefficients $\frac{r_is_i-1}{r_is_i}$ to the irreducible components $E_i$,
where the $s_i$ are defined in the proof of the above lemma.
We note that, even if we start with the non-boundary case $B = 0$, 
the naturally defined boundary divisor $C$ on $Y$ is non-zero in general, 
because there may be multiple fibers for $\phi$.
Let $\mathcal{Y}$ be the smooth Deligne-Mumford stack 
associated to the pair $(Y,C)$.
The above lemma implies the following:

\begin{Cor}\label{smooth}
The natural morphism $\psi: \mathcal{X} \to \mathcal{Y}$ is smooth.
\end{Cor}

\begin{Thm}\label{fiber}
(1) The functor $\psi^*: D^b(\text{Coh}(\mathcal{Y})) \to 
D^b(\text{Coh}(\mathcal{X}))$ is fully faithful. 

Let $D^b(\text{Coh}(\mathcal{Y}))_k$ denote the full subcategory of 
$D^b(\text{Coh}(\mathcal{X}))$ defined by
\[
D^b(\text{Coh}(\mathcal{Y}))_k = 
\psi^*D^b(\text{Coh}(\mathcal{Y})) \otimes 
\mathcal{O}_{\mathcal{X}}(\sum_{i=1}^{\alpha} k_i \mathcal{D}_i)
\]
for a sequence of integers $k = (k_i)$ for $1 \le i \le \alpha$.

(2) If $\sum_{i=1}^{\alpha} \frac{a_ik_i}{r_i} > 
\sum_{i=1}^{\alpha} \frac{a_ik'_i}{r_i} > 
\sum_{i=1}^{\alpha} \frac{a_i(k_i-1)}{r_i}$, then
\[
\text{Hom}^q(D^b(\text{Coh}(\mathcal{Y}))_k,
D^b(\text{Coh}(\mathcal{Y}))_{k'}) = 0
\]
for all $q$, where $k' = (k'_i)$ is another sequence of integers.

(3) If $\sum_{i=1}^{\alpha} \frac{a_ik_i}{r_i} = 
\sum_{i=1}^{\alpha} \frac{a_ik'_i}{r_i}$ and
$\mathcal{O}_{\mathcal{X}}(\sum_{i=1}^{\alpha} (k_i-k'_i) \mathcal{D}_i)
\not\in \psi^*D^b(\text{Coh}(\mathcal{Y}))$, then
\[
\text{Hom}^q(D^b(\text{Coh}(\mathcal{Y}))_k,
D^b(\text{Coh}(\mathcal{Y}))_{k'}) = 0
\]
for all $q$.

(4) The set of subcategories
$D^b(\text{Coh}(\mathcal{Y}))_k$ for 
\[
0 \ge \sum_{i=1}^{\alpha} \frac{a_ik_i}{r_i} > 
- \sum_{i=1}^{\alpha} \frac{a_i}{r_i}
\]
generates the triangulated category $D^b(\text{Coh}(\mathcal{X}))$.
\end{Thm}

\begin{proof}
(1) By \cite{BO} or \cite{Bridgeland1}, it is sufficient to prove the following
statement; if $A$ and $B$ are skyscraper sheaves on $\mathcal{Y}$ of 
length $1$, then the natural homomorphism
$\text{Hom}^p(A,B) \to \text{Hom}^p(\psi^*A,\psi^*B)$ is bijective.
This follows from the fact that 
$X'_0$ is isomorphic to the product of $Y''_0$ 
with a weighted projective space $\mathbb{P}(a_1,\dots,a_{\alpha})$, and that
the natural homomorphism of Galois groups
$N_X/N_{X'_0} \to N_Y/N_{Y''_0}$ is surjective .

For (2) and (3), we use a spectral sequence
\[
\begin{split}
&E_2^{p,q} = H^p(\mathcal{Y}, \mathcal{H}om(A,B) \otimes 
R^q\psi_* \mathcal{O}_{\mathcal{X}}(\sum_{i=1}^{\alpha} k_i \mathcal{D}_i)) \\
&\Rightarrow \text{Hom}^{p+q}(\psi^*A,\psi^*B 
\otimes \mathcal{O}_{\mathcal{X}}(\sum_{i=1}^{\alpha} k_i \mathcal{D}_i)) 
\end{split}
\]
for invertible sheaves $A,B$ on $\mathcal{Y}$.
The direct image sheaves vanish in our case, 
because the relative canonical divisor for $\psi$ is given by 
\[
\omega_{\mathcal{X}/\mathcal{Y}} \cong \mathcal{O}_{\mathcal{X}}
(- \sum_{i=1}^{\alpha} \mathcal{D}_i)
\]
and an invertible sheaf 
$\mathcal{O}_{\mathcal{X}}(\sum_{i=1}^{\alpha} k_i \mathcal{D}_i)$
is $\psi$-ample if and only if $\sum_{i=1}^{\alpha} \frac{a_ik_i}{r_i} > 0$
(\cite{KMM}).

(4) In general, a full triangulated subcategory $\mathcal{B}$ 
of a triangulated category $\mathcal{A}$ is said to be 
{\em right} (resp. {\em left}) {\em admissible} if $\mathcal{A}$ is generated 
by $\mathcal{B}$ and $\mathcal{B}^{\perp}$ (resp. $\mathcal{B}$ 
and ${}^{\perp}\mathcal{B}$), where $\mathcal{B}^{\perp}$ 
(resp. ${}^{\perp}\mathcal{B}$) denotes the right (resp. left) orthogonal 
complement of $\mathcal{B}$ in $\mathcal{A}$ (\cite{BK}).
The subcategory triangulated $\mathcal{T}$ of $D^b(\text{Coh}(\mathcal{X}))$ 
generated by the subcategories $D^b(\text{Coh}(\mathcal{Y}))_k$ 
is admissible by 
loc.~cit.~1.12, 2.6 and 2.11.
Therefore, it is sufficient to prove that the left orthogonal 
${}^{\perp}\mathcal{T}$ consists of $0$ objects.

Let $A$ be an arbitrary skyscraper sheaf of length $1$ on $\mathcal{X}$
supported at a point $P$.
Then by Theorem~\ref{Fano}, there exists a skyscraper sheaf $B$ of length 
$1$ on $\mathcal{Y}$ supported at $Q = \phi(P)$ 
such that $A$ is contained in the subcategory 
generated by the sheaves of the form 
$\psi^*B \otimes \mathcal{O}_{\mathcal{X}}
(\sum_{i=1}^{\alpha} k_i\mathcal{D}_i)$ for 
$0 \ge \sum_{i=1}^{\alpha} \frac{a_ik_i}{r_i} > 
- \sum_{i=1}^{\alpha} \frac{a_i}{r_i}$.
Thus $A$ contained in $\mathcal{T}$.
Hence ${}^{\perp}\mathcal{T} = 0$,
because such $A$ span $D^b(\text{Coh}(\mathcal{X}))$ 
(\cite{BO} or \cite{Bridgeland1}).
\end{proof}

\begin{Cor}\label{cfiber}
Assume that $D^b(\text{Coh}(\mathcal{Y}))$ has a complete exceptional 
collection consisting of sheaves. 
Then so has $D^b(\text{Coh}(\mathcal{X}))$.
\end{Cor}


\section{Divisorial contraction}

We consider a toric divisorial contraction $\phi: X \to Y$.
$K_X+B$ is negative for $\phi$, and $C = \phi_*B$ is the strict transform.
Let $D$ be the exceptional divisor of the contraction.
Then the restriction $\bar{\phi}: D \to F = \phi(D)$ is a Mori fiber space 
which was treated in the previous section.

Let $\mathcal{Y}$ be the stack associated to the pair $(Y,C)$.
We note that there is no morphism of stacks from 
$\mathcal{X}$ to $\mathcal{Y}$ in general.
But there is still a fully faithful functor
$\Phi: D^b(\text{Coh}(\mathcal{Y})) \to D^b(\text{Coh}(\mathcal{X}))$ by
\cite{log-crep}~Theorem~4.2~(2).
Indeed, let $\mathcal{W}$ be the normalization of the fiber product 
$\mathcal{X} \times_Y \mathcal{Y}$, and let $\mu: \mathcal{W} \to \mathcal{X}$
and $\nu: \mathcal{W} \to \mathcal{Y}$ be the projections.
Then $\Phi = \mu_* \circ \nu^*$ is fully faithful.
We regard $D^b(\text{Coh}(\mathcal{Y}))$ as a full subcategory of 
$D^b(\text{Coh}(\mathcal{X}))$ through this functor.

Let $E_i = \phi_*D_i$ be the prime divisors on $Y$ corresponding 
to the edges $v_i$ for $1 \le i \le n$.
Those $E_i$ for $1 \le i \le \alpha$ are the divisors 
which contain the center $F$ of the blowing-up $\phi$,
and $D = D_{n+1}$ is the exceptional divisor.
Let $\mathcal{E}_i$ be the prime divisors on $\mathcal{Y}$ corresponding to 
the $E_i$.
The following formula is proved in the proof of 
\cite{log-crep}~Theorem~4.2~(2):
\[
\begin{split}
&\Phi(\mathcal{O}_{\mathcal{Y}}(\sum_{i=1}^n k_i\mathcal{E}_i))
\cong \mathcal{O}_{\mathcal{X}}(\sum_{i=1}^{n+1} k_i\mathcal{D}_i) \\
&k_{n+1} = \llcorner \frac{r_{n+1}}{b_{n+1}} \sum_{i=1}^n \frac{a_ik_i}{r_i} 
\lrcorner
\end{split}
\]
for any integers $k_i$ for $1 \le i \le n$, 
where we put $b_{n+1} = - a_{n+1} > 0$.

Let $r$ be a positive integer such that $a_ir$ is divisible by $r_i$ for 
$1 \le i \le n+1$.
We set 
\[
\vert a_i \vert \cdot r = r_is_i.
\]
Let $s = (s_1, \dots, s_{n+1})$ be the 
greatest common divisor, and set $s_i = s \bar s_i$.
Then the fractional part of the rational number
\[
\frac{r_{n+1}}{b_{n+1}} \sum_{i=1}^n \frac{a_ik_i}{r_i}
= \frac{\sum_{i=1}^n k_i \bar s_i}{\bar s_{n+1}}
\]
can take arbitrary value in the set 
$\{0, \frac 1{\bar s_{n+1}}, \dots, \frac{\bar s_{n+1}-1}{\bar s_{n+1}}\}$
when we vary the sequence $k$,
because $(\bar s_1, \dots, \bar s_{n+1}) = 1$.

The Mori fiber space $\bar{\phi}: D \to F$ is described as follows.
The lattice of $1$-parameter subgroup for $D$ is given by 
$\bar{N} = N_X/\mathbb{Z}v_{n+1}$.
We write $v_i \text{ mod }\mathbb{Z}v_{n+1} = t_i \bar v_i$
for $1 \le i \le n$, where the 
$t_i$ are positive integers and the $\bar v_i$ are primitive vectors 
in $\bar{N}$.
Let $t = (a_1t_1, \dots, a_nt_n)$ be the greatest common divisor, 
and denote $a_it_i = t\bar a_i$.
Then we have an equation
\[
\bar a_1\bar v_1 + \dots + \bar a_n\bar v_n = 0.
\]

We define a $\mathbb{Q}$-divisor $\bar B$ on $D$ by putting coefficients
$\frac{r_it_i-1}{r_it_i}$ to the prime divisors $\bar D_i = D_i \cap D$
for $1 \le i \le n$. 
We also define a $\mathbb{Q}$-divisor $\bar C$ on the base space of 
the Mori fiber space $F$ using $\bar B$ as in the previous section.
Let $\mathcal{D}$ and $\mathcal{F}$ be the smooth stacks associated to the 
pairs $(D,\bar B)$ and $(F,\bar C)$, respectively.
Then there are induced morphisms of stacks 
$\bar{\psi}: \mathcal{D} \to \mathcal{F}$ and $j: \mathcal{D} \to \mathcal{X}$.
Let $\bar{\mathcal{D}}_i$ be the prime divisors on $\mathcal{D}$ 
corresponding to the $\bar D_i$ for $1 \le i \le n$. 
Then we have 
\[
j^*\mathcal{O}_{\mathcal{X}}(\mathcal{D}_i)
\cong \mathcal{O}_{\mathcal{D}}(\bar{\mathcal{D}}_i).
\]
We note that $D_i \vert_D = \frac 1{t_i} \bar D_i$ in the usual language.

\begin{Rem}\label{stabilizer}
If $r_{n+1} > 1$, then the action of the stabilizer group at the generic point
of $\mathcal{D}_{n+1}$ is non-trivial.
Hence we have $j^*\mathcal{O}_{\mathcal{X}}(k\mathcal{D}_{n+1}) = 0$ 
on $\mathcal{D}$ if $k$ is not divisible by $r_{n+1}$.
Indeed, we have 
\[
\text{Hom}(j^*\mathcal{O}_{\mathcal{X}}(k\mathcal{D}_{n+1}), A)
\cong \text{Hom}(\mathcal{O}_{\mathcal{X}}(k\mathcal{D}_{n+1}), j_*A)
\cong 0
\]
for any sheaf $A$ on $\mathcal{D}$ in this case.

For example, let $X$ be an affine line with a point $P$, and $\mathcal{X}$ 
the stack associated to the pair $(X, \frac{r-1}r P)$ with a point 
$\mathcal{P}$ above $P$.
Then we have $j^*\mathcal{O}_{\mathcal{X}}(k\mathcal{P}) = 0$ 
if $k$ is not divisible by $r$, 
where $j: P \to \mathcal{X}$ is the natural morphism.
From a resolution
\[
0 \to \mathcal{O}_{\mathcal{X}}((k-1)\mathcal{P}) \to 
\mathcal{O}_{\mathcal{X}}(k\mathcal{P}) 
\to \mathcal{O}_{\mathcal{P}}(k\mathcal{P}) \to 0
\]
it follows that $L_qj^*\mathcal{O}_{\mathcal{P}}(k\mathcal{P})$ is isomorphic 
to $\mathcal{O}_P$ if $q=0$ and $k \equiv 0  \mod r$, or 
$q = 1$ and $k \equiv 1 \mod r$, and $0$ if otherwise.
Thus
\[
\begin{split}
&\text{Hom}^q(\mathcal{O}_{\mathcal{P}}(k\mathcal{P}), 
\mathcal{O}_{\mathcal{P}})
\cong \text{Hom}^q(\mathcal{O}_{\mathcal{P}}(k\mathcal{P}), j_*\mathcal{O}_P) 
\\
&\cong \text{Hom}^q(Lj^*\mathcal{O}_{\mathcal{P}}(k\mathcal{P}), \mathcal{O}_P)
\end{split}
\]
is non zero if and only if $q=0$ and $k \equiv 0 \mod r$, or 
$q = 1$ and $k \equiv 1 \mod r$.
\end{Rem}

\begin{Thm}\label{divisorial}
(1) The functor $j_*\bar{\psi}^*: D^b(\text{Coh}(\mathcal{F})) \to 
D^b(\text{Coh}(\mathcal{X}))$ is fully faithful. 

Let $D^b(\text{Coh}(\mathcal{F}))_k$ denote the full subcategory of 
$D^b(\text{Coh}(\mathcal{X}))$ defined by 
\[
D^b(\text{Coh}(\mathcal{F}))_k = j_*\bar{\psi}^*D^b(\text{Coh}(\mathcal{F})) 
\otimes \mathcal{O}_{\mathcal{X}}(\sum_{i=1}^{n+1} k_i \mathcal{D}_i)
\]
for a sequence of integers $k = (k_i)$ for $1 \le i \le n+1$. 

(2) If $0 > \sum_{i=1}^{n+1} \frac{a_ik_i}{r_i} 
\ge - \sum_{i=1}^{n+1} \frac{a_i}{r_i}$, then
\[
\text{Hom}^q(\Phi(D^b(\text{Coh}(\mathcal{Y}))),
D^b(\text{Coh}(\mathcal{F}))_k) = 0
\]
for all $q$.

(3) If $\sum_{i=1}^{n+1} \frac{a_ik_i}{r_i} > 
\sum_{i=1}^{n+1} \frac{a_ik'_i}{r_i} > 
\sum_{i=1}^{n+1} \frac{a_i(k_i-1)}{r_i}$, then
\[
\text{Hom}^q(D^b(\text{Coh}(\mathcal{F}))_k,
D^b(\text{Coh}(\mathcal{F}))_{k'}) = 0
\]
for all $q$, where $k' = (k'_i)$ is another sequence of integers.

(4) If $\sum_{i=1}^{n+1} \frac{a_ik_i}{r_i} = 
\sum_{i=1}^{n+1} \frac{a_ik'_i}{r_i}$, but if
$j^*\mathcal{O}_{\mathcal{X}}(\sum_{i=1}^{n+1} (k_i-k'_i) \mathcal{D}_i)=0$
or 
\[
j^*\mathcal{O}_{\mathcal{X}}(\sum_{i=1}^{n+1} (k_i-k'_i) \mathcal{D}_i)
\not\in \bar{\psi}^*D^b(\text{Coh}(\mathcal{F}))
\]
then
\[
\text{Hom}^q(D^b(\text{Coh}(\mathcal{F}))_k,
D^b(\text{Coh}(\mathcal{F}))_{k'}) = 0
\]
for all $q$.

(5) The subcategories $\Phi(D^b(\text{Coh}(\mathcal{Y})))$ and the 
$D^b(\text{Coh}(\mathcal{F}))_k$ for 
\[
0 > \sum_{i=1}^{n+1} \frac{a_ik_i}{r_i} \ge - \sum_{i=1}^{n+1} \frac{a_i}{r_i}
\]
generate the triangulated category $D^b(\text{Coh}(\mathcal{X}))$.
\end{Thm}

\begin{proof}
(1) It is sufficient to prove that the natural homomorphism
\[
\text{Hom}^q(L,L') \to \text{Hom}^q(j_*\bar{\psi}^*L,j_*\bar{\psi}^*L')
\]
is bijective for all $q$ and all locally free sheaves $L$ and $L'$ 
on $\mathcal{F}$, because these sheaves span the category 
$D^b(\text{Coh}(\mathcal{F}))$.

We have an exact sequence 
\[
0 \to \mathcal{O}_{\mathcal{X}}(- \mathcal{D}_{n+1}) \to 
\mathcal{O}_{\mathcal{X}} \to \mathcal{O}_{\mathcal{D}_{n+1}} \to 0
\]
with an isomorphism $\mathcal{O}_{\mathcal{D}_{n+1}}
\cong j_*\mathcal{O}_{\mathcal{D}}$. 
Hence 
\[
L_qj^*j_*\mathcal{O}_{\mathcal{D}}
\cong \begin{cases} \mathcal{O}_{\mathcal{D}} &\text{ for } q = 0 \\ 
j^*\mathcal{O}_{\mathcal{X}}(- \mathcal{D}_{n+1}) &\text{ for } q = 1 \\
0 &\text{ otherwise}
\end{cases}
\]
where $j^*\mathcal{O}_{\mathcal{X}}(- \mathcal{D}_{n+1})$ is an invertible 
sheaf on $\mathcal{D}$ if $r_{n+1} = 1$, and zero otherwise.

If $r_{n+1} > 1$, then
\[
\begin{split}
&\text{Hom}^q(j_*\bar{\psi}^*L,j_*\bar{\psi}^*L')
\cong \text{Hom}^q(Lj^*j_*\bar{\psi}^*L,\bar{\psi}^*L') \\
&\cong \text{Hom}^q(\bar{\psi}^*L,\bar{\psi}^*L')
\cong \text{Hom}^q(L,L')
\end{split}
\]
as required.
If $r_{n+1} = 1$, then we know that 
$j^*\mathcal{O}_{\mathcal{X}}(\mathcal{D}_{n+1})$ 
is negative for $\bar{\psi}$, while 
$j^*\mathcal{O}_{\mathcal{X}}(\sum_{i=1}^{n+1} \mathcal{D}_{n+1})$ is ample 
for $\bar{\psi}$, because $\sum_{i=1}^{n+1} \frac{a_i}{r_i} > 0$. 
Since 
\[
\omega_{\mathcal{D}/\mathcal{F}} 
\cong \mathcal{O}_{\mathcal{D}}(- \sum_{i=1}^{\alpha} \bar{\mathcal{D}}_i)
\cong j^*\mathcal{O}_{\mathcal{X}}(- \sum_{i=1}^{\alpha} \mathcal{D}_i)
\]
we calculate
\[
\text{Hom}^q(L_1j^*j_*\bar{\psi}^*L,\bar{\psi}^*L')
\cong \text{Hom}^q(\bar{\psi}^*L,
\bar{\psi}^*L' \otimes j^*\mathcal{O}_{\mathcal{X}}(\mathcal{D}_{n+1}))
\cong 0
\]
by the relative vanishing theorem for $\bar{\psi}$ (\cite{KMM}).
Therefore, we have also our assertion in this case.

(2) It is sufficient to prove 
\[
\text{Hom}^q(\mathcal{O}_{\mathcal{X}}(\sum_{i=1}^{n+1} k'_i\mathcal{D}_i),
j_*\bar{\psi}^*A \otimes 
\mathcal{O}_{\mathcal{X}}(\sum_{i=1}^{n+1} k_i \mathcal{D}_i))
= 0
\]
for all integers $q$, all sheaves $A$ on $\mathcal{F}$, 
and for the sequences $(k)$ and $(k')$ under the additional conditions that 
\[
\begin{split}
&k'_{n+1} = \llcorner \frac{r_{n+1}}{b_{n+1}} \sum_{i=1}^n \frac{a_ik'_i}{r_i} 
\lrcorner \\
&0 > \sum_{i=1}^{n+1} \frac{a_ik_i}{r_i} 
\ge - \sum_{i=1}^{n+1} \frac{a_i}{r_i}.
\end{split}
\]
By the first condition, we have
\[
0 \le \sum_{i=1}^{n+1} \frac{a_i k'_i}{r_i} < \frac {b_{n+1}}{r_{n+1}}.
\]
Hence
\[
0 > \sum_{i=1}^{n+1} \frac{a_i(k_i - k'_i)}{r_i} 
> - \sum_{i=1}^n \frac{a_i}{r_i}.
\]
By the relative vanishing theorem for $\bar{\psi}$, we have 
\[
\begin{split}
&\text{Hom}^q(\mathcal{O}_{\mathcal{X}}(\sum_{i=1}^{n+1} k'_i\mathcal{D}_i),
j_*\bar{\psi}^*A \otimes 
\mathcal{O}_{\mathcal{X}}(\sum_{i=1}^{n+1} k_i \mathcal{D}_i)) \\
&\cong \text{Hom}^q(j^*\mathcal{O}_{\mathcal{X}}
(\sum_{i=1}^{n+1} (k'_i-k_i)\mathcal{D}_i), \bar{\psi}^*A)
\cong 0.
\end{split}
\]

(3) is similarly proved as in (1).
Since $0 > \sum_{i=1}^{n+1} \frac{a_i(k'_i-k_i)}{r_i} > 
- \sum_{i=1}^n \frac{a_i}{r_i} + \frac{b_{n+1}}{r_{n+1}}$, we have
\[
R\bar{\psi}_*j^*\mathcal{O}_{\mathcal{X}}
(\sum_{i=1}^{n+1} (k'_i-k_i)\mathcal{D}_i)
= R\bar{\psi}_*j^*\mathcal{O}_{\mathcal{X}}
(\mathcal{D}_{n+1} + \sum_{i=1}^{n+1} (k'_i-k_i)\mathcal{D}_i)
= 0
\]
by the relative vanishing theorem for $\bar{\psi}$.
Thus
\[
\text{Hom}^q(j_*\bar{\psi}^*L,j_*\bar{\psi}^*L' \otimes 
\mathcal{O}_{\mathcal{X}}
(\sum_{i=1}^{n+1} (k'_i-k_i)\mathcal{D}_i)) = 0
\]
for all $q$ and all locally free sheaves $L$ and $L'$ 
on $\mathcal{F}$.

(4) is similar to (3).

(5) We shall prove that the left orthogonal ${}^{\perp}\mathcal{T}$ to 
the triangulated subcategory $\mathcal{T}$ of $D^b(\text{Coh}(\mathcal{X}))$ 
generated by these subcategories consists of $0$ objects as in the proof 
of Theorem~\ref{fiber}.

Let $A$ be an arbitrary skyscraper sheaf of length $1$ on $\mathcal{X}$
supported at a point $P$.
If $P \not\in \mathcal{D}_{n+1}$, then $A \in \mathcal{T}$.
Otherwise, there is a point $\bar P$ on $\mathcal{D}$ such that 
$P = j(\bar P)$.
Then by Theorem~\ref{Fano}, there exists a skyscraper sheaf $B$ of length 
$1$ on $\mathcal{F}$ supported at $\bar Q = \bar{\psi}(\bar P)$ 
such that $A$ is contained in the subcategory 
generated by the sheaves of the form 
$j_*\bar{\psi}^*B \otimes \mathcal{O}_{\mathcal{X}}
(\sum_{i=1}^{n+1} k_i\mathcal{D}_i)$
for 
\[
\frac{b_{n+1}}{r_{n+1}} > \sum_{i=1}^{n+1} \frac{a_ik_i}{r_i} \ge 
- \sum_{i=1}^{n+1} \frac{a_i}{r_i}.
\]
If $\frac{b_{n+1}}{r_{n+1}} > \sum_{i=1}^{n+1} \frac{a_ik_i}{r_i} \ge 0$, then
it follows that 
$k_{n+1} = \llcorner \frac{r_{n+1}}{b_{n+1}} \sum_{i=1}^n \frac{a_ik_i}{r_i} 
\lrcorner$.
Therefore, $A$ contained in $\mathcal{T}$, hence ${}^{\perp}\mathcal{T} = 0$.
\end{proof}

\begin{Cor}\label{cdivisorial}
Assume that $D^b(\text{Coh}(\mathcal{Y}))$ has a complete exceptional 
collection consisting of sheaves. 
Then so has $D^b(\text{Coh}(\mathcal{X}))$.
\end{Cor}


\section{Log flip}

We consider a toric small contraction $\phi: X \to Y$ with the
log flip $\phi^+: X^+ \to Y$.
$K_X+B$ is negative for $\phi$, and $K_{X^+}+B^+$ is ample for $\phi^+$, 
where $B^+ = (\phi^+_*)^{-1}\phi_*B$ is the strict transform.
The argument for log flips in this section is surprizingly similar to 
that for the divisorial contractions in the previous section.

Let $\mathcal{X}^+$ be the smooth Deligne-Mumford stack 
associated to the pair $(X^+,B^+)$.
Then there is a fully faithful functor
$\Phi: D^b(\text{Coh}(\mathcal{X}^+)) \to D^b(\text{Coh}(\mathcal{X}))$ by
\cite{log-crep}~Theorem~4.2~(3).
Indeed, let $\mathcal{W}$ be the normalization of the fiber product 
$\mathcal{X} \times_Y \mathcal{X}^+$, 
and let $\mu: \mathcal{W} \to \mathcal{X}$
and $\nu: \mathcal{W} \to \mathcal{X}^+$ be the projections.
Then $\Phi = \mu_* \circ \nu^*$ is fully faithful.
We regard $D^b(\text{Coh}(\mathcal{X}^+))$ as a full subcategory of 
$D^b(\text{Coh}(\mathcal{X}))$ through this functor.

Let $D^+_i = (\phi^+_*)^{-1}\phi_*D_i$ be the prime divisors on $X^+$ 
corresponding to the edges $v_i$ for $1 \le i \le n+1$, and 
let $\mathcal{D}^+_i$ be the corresponding prime divisors on $\mathcal{X}^+$.
The following formula is proved in the proof of 
\cite{log-crep}~Theorem~4.2~(3):
\[
\Phi(\mathcal{O}_{\mathcal{X}^+}(\sum_{i=1}^{n+1} k_i\mathcal{D}^+_i)) 
\cong \mathcal{O}_{\mathcal{X}}(\sum_{i=1}^{n+1} k_i\mathcal{D}_i)
\]
if
\[
0 \le \sum_{i=1}^{n+1} \frac{a_ik_i}{r_i} < 
\sum_{i=\beta+1}^{n+1} \frac{b_i}{r_i}
\]
where we put $b_i = - a_i$ for $\beta+1 \le i \le n+1$.

Let $D$ be the exceptional locus of the contraction $\phi$.
Then we have $D = \cap_{i=\beta+1}^{n+1} D_i$, and 
the restriction $\bar{\phi}: D \to F = \phi(D)$ is a Mori fiber space,
which is described as follows.
The lattice of $1$-parameter subgroup for $D$ is given by 
$\bar{N} = N_X/\bigoplus_{i=\beta+1}^{n+1} \mathbb{Z}v_i$.
We write $v_i \text{ mod }\bigoplus_{i=\beta+1}^{n+1} \mathbb{Z}v_i 
= t_i \bar v_i$ for $1 \le i \le \beta$, where the
$t_i$ are positive integers and the $\bar v_i$ are primitive vectors 
in $\bar{N}$.
Let $t = (a_1t_1, \dots, a_{\beta}t_{\beta})$ be the greatest common divisor
and $a_it_i = t\bar a_i$.
Then we have an equation
\[
\bar a_1\bar v_1 + \dots + \bar a_{\beta}\bar v_{\beta} = 0.
\]

We define a $\mathbb{Q}$-divisor $\bar B$ on $D$ by putting coefficients
$\frac{r_it_i-1}{r_it_i}$ to the prime divisors $\bar D_i = D_i \cap D$
for $1 \le i \le \beta$. 
We also define a $\mathbb{Q}$-divisor $\bar C$ on the base space of 
the Mori fiber space $F$ using $\bar B$ as before.
Let $\mathcal{D}$ and $\mathcal{F}$ be the smooth Deligne-Mumford 
stacks associated to the 
pairs $(D,\bar B)$ and $(F,\bar C)$, respectively.
Then there are induced morphisms of stacks 
$\bar{\psi}: \mathcal{D} \to \mathcal{F}$ and $j: \mathcal{D} \to \mathcal{X}$.
Let $\bar{\mathcal{D}}_i$ be the prime divisors on $\mathcal{D}$ 
corresponding to the $\bar D_i$ for $1 \le i \le \beta$. 
Then we have 
\[
j^*\mathcal{O}_{\mathcal{X}}(\mathcal{D}_i)
\cong \mathcal{O}_{\mathcal{D}}(\bar{\mathcal{D}}_i).
\]
We note that $D_i \vert_D = \frac 1{t_i} \bar D_i$ in the usual language.

\begin{Thm}\label{flip}
(1) The functor $j_*\bar{\psi}^*: D^b(\text{Coh}(\mathcal{F})) \to 
D^b(\text{Coh}(\mathcal{X}))$ is fully faithful. 

Let $D^b(\text{Coh}(\mathcal{F}))_k$ denote the full subcategory of 
$D^b(\text{Coh}(\mathcal{X}))$ defined by 
\[
D^b(\text{Coh}(\mathcal{F}))_k = j_*\bar{\psi}^*D^b(\text{Coh}(\mathcal{F})) 
\otimes \mathcal{O}_{\mathcal{X}}(\sum_{i=1}^{n+1} k_i \mathcal{D}_i)
\]
for a sequence of integers $k = (k_i)$ for $1 \le i \le n+1$. 

(2) If $0 > \sum_{i=1}^{n+1} \frac{a_ik_i}{r_i} 
\ge - \sum_{i=1}^{n+1} \frac{a_i}{r_i}$, then
\[
\text{Hom}^q(\Phi(D^b(\text{Coh}(\mathcal{X}^+))),
D^b(\text{Coh}(\mathcal{F}))_k) = 0
\]
for all $q$.

(3) If $\sum_{i=1}^{n+1} \frac{a_ik_i}{r_i} > 
\sum_{i=1}^{n+1} \frac{a_ik'_i}{r_i} > 
\sum_{i=1}^{n+1} \frac{a_i(k_i-1)}{r_i}$, then
\[
\text{Hom}^q(D^b(\text{Coh}(\mathcal{F}))_k,
D^b(\text{Coh}(\mathcal{F}))_{k'}) = 0
\]
for all $q$, where $k' = (k'_i)$ is another sequence of integers.

(4) If $\sum_{i=1}^{n+1} \frac{a_ik_i}{r_i} = 
\sum_{i=1}^{n+1} \frac{a_ik'_i}{r_i}$, but if
$j^*\mathcal{O}_{\mathcal{X}}(\sum_{i=1}^{n+1} (k_i-k'_i) \mathcal{D}_i)=0$
or 
\[
j^*\mathcal{O}_{\mathcal{X}}(\sum_{i=1}^{n+1} (k_i-k'_i) \mathcal{D}_i)
\not\in \bar{\psi}^*D^b(\text{Coh}(\mathcal{F}))
\]
then
\[
\text{Hom}^q(D^b(\text{Coh}(\mathcal{F}))_k,
D^b(\text{Coh}(\mathcal{F}))_{k'}) = 0
\]
for all $q$.

(5) The subcategories $\Phi(D^b(\text{Coh}(\mathcal{X}^+)))$ and the 
$D^b(\text{Coh}(\mathcal{F}))_k$ for 
\[
0 > \sum_{i=1}^{n+1} \frac{a_ik_i}{r_i} \ge - \sum_{i=1}^{n+1} \frac{a_i}{r_i}
\]
generate the triangulated category $D^b(\text{Coh}(\mathcal{X}))$.
\end{Thm}

\begin{proof}
(1) We shall prove that the natural homomorphism
\[
\text{Hom}^q(L,L') \to \text{Hom}^q(j_*\bar{\psi}^*L,j_*\bar{\psi}^*L')
\]
is bijective for all $q$ and all locally free sheaves $L$ and $L'$ 
on $\mathcal{F}$.

We have an exact sequence 
\[
\begin{split}
&0 \to \mathcal{O}_{\mathcal{X}}(- \sum_{i=\beta+1}^{n+1} \mathcal{D}_i) 
\to \cdots \to 
\bigwedge^2 (\bigoplus_{i=\beta+1}^{n+1} \mathcal{O}_{\mathcal{X}}
(- \mathcal{D}_i)) \\
&\to \bigoplus_{i=\beta+1}^{n+1} \mathcal{O}_{\mathcal{X}}(- \mathcal{D}_i) 
\to \mathcal{O}_{\mathcal{X}} \to j_*\mathcal{O}_{\mathcal{D}} \to 0.
\end{split}
\]
Hence 
\[
L_qj^*j_*\mathcal{O}_{\mathcal{D}}
\cong \bigwedge^q (\bigoplus_{i=\beta+1}^{n+1} j^*\mathcal{O}_{\mathcal{X}}
(- \mathcal{D}_i)).
\]
The sheaf $j^*\mathcal{O}_{\mathcal{X}}(- \sum_{i \in I} \mathcal{D}_i)$ 
for any subset $I \subset \{\beta+1, \dots, n+1\}$ is either invertible 
or zero, and is negative for $\bar{\psi}$ if it is not a zero sheaf. 

Since 
\[
\omega_{\mathcal{D}/\mathcal{F}} 
\cong \mathcal{O}_{\mathcal{D}}(- \sum_{i=1}^{\alpha} \bar{\mathcal{D}}_i)
\cong j^*\mathcal{O}_{\mathcal{X}}(- \sum_{i=1}^{\alpha} \mathcal{D}_i)
\]
we calculate
\[
\text{Hom}^q(L_pj^*j_*\bar{\psi}^*L,\bar{\psi}^*L')
\cong \text{Hom}^q(\bar{\psi}^*L,
\bar{\psi}^*L' \otimes 
\bigwedge^p (\bigoplus_{i=\beta+1}^{n+1} j^*\mathcal{O}_{\mathcal{X}}
(- \mathcal{D}_i)))
\cong 0
\]
for $p > 0$ and for any $q$ by the relative vanishing theorem for $\bar{\psi}$,
because $\sum_{i=1}^{n+1} \frac{a_i}{r_i} > 0$. 
Hence 
\[
\begin{split}
&\text{Hom}^q(j_*\bar{\psi}^*L,j_*\bar{\psi}^*L')
\cong \text{Hom}^q(j^*j_*\bar{\psi}^*L,\bar{\psi}^*L') \\
&\cong \text{Hom}^q(\bar{\psi}^*L,\bar{\psi}^*L')
\cong \text{Hom}^q(L,L')
\end{split}
\]
for any $q$ as required.

(2) It is sufficient to prove 
\[
\text{Hom}^q(\mathcal{O}_{\mathcal{X}}(\sum_{i=1}^{n+1} k'_i\mathcal{D}_i),
j_*\bar{\psi}^*A \otimes 
\mathcal{O}_{\mathcal{X}}(\sum_{i=1}^{n+1} k_i \mathcal{D}_i))
= 0
\]
for all integers $q$, all sheaves $A$ on $\mathcal{F}$, 
and for the sequences $(k)$ and $(k')$ under the additional conditions that 
\[
\begin{split}
&0 \le \sum_{i=1}^{n+1} \frac{a_ik'_i}{r_i} 
< \sum_{i=\beta+1}^{n+1} \frac{b_i}{r_i} \\
&0 > \sum_{i=1}^{n+1} \frac{a_ik_i}{r_i} 
\ge - \sum_{i=1}^{n+1} \frac{a_i}{r_i}.
\end{split}
\]
It follows that
\[
0 > \sum_{i=1}^{n+1} \frac{a_i(k_i - k'_i)}{r_i} 
> - \sum_{i=1}^{\alpha} \frac{a_i}{r_i}.
\]
By the relative vanishing theorem for $\bar{\psi}$, we have 
\[
\begin{split}
&\text{Hom}^q(\mathcal{O}_{\mathcal{X}}(\sum_{i=1}^{n+1} k'_i\mathcal{D}_i),
j_*\bar{\psi}^*A \otimes 
\mathcal{O}_{\mathcal{X}}(\sum_{i=1}^{n+1} k_i \mathcal{D}_i)) \\
&\cong \text{Hom}^q(j^*\mathcal{O}_{\mathcal{X}}
(\sum_{i=1}^{n+1} (k'_i-k_i)\mathcal{D}_i), \bar{\psi}^*A)
\cong 0.
\end{split}
\]

(3) is similarly proved as in (1).
Since $0 > \sum_{i=1}^{n+1} \frac{a_i(k'_i-k_i)}{r_i} > 
- \sum_{i=1}^{\alpha} \frac{a_i}{r_i} + 
\sum_{i=\beta+1}^{n+1} \frac{b_i}{r_i}$, we have
\[
R\bar{\psi}_*j^*\mathcal{O}_{\mathcal{X}}
(\sum_{i \in I} \mathcal{D}_i + \sum_{i=1}^{n+1} (k'_i-k_i)\mathcal{D}_i)
= 0
\]
for any subset $I \subset \{\beta+1, \dots, n+1\}$ 
by the relative vanishing theorem for $\bar{\psi}$.
Thus
\[
\text{Hom}^q(j_*\bar{\psi}^*L,j_*\bar{\psi}^*L' \otimes 
\mathcal{O}_{\mathcal{X}}
(\sum_{i=1}^{n+1} (k'_i-k_i)\mathcal{D}_i)) = 0
\]
for all $q$ and all locally free sheaves $L$ and $L'$ 
on $\mathcal{F}$.

(4) is similar to (3).

(5) We shall prove that the left orthogonal ${}^{\perp}\mathcal{T}$ to 
the triangulated subcategory $\mathcal{T}$ of $D^b(\text{Coh}(\mathcal{X}))$ 
generated by these subcategories consists of $0$ objects.

Let $A$ be an arbitrary skyscraper sheaf of length $1$ on $\mathcal{X}$
supported at a point $P$.
If $P$ is not above a point in $D$, then $A \in \mathcal{T}$.
Otherwise, there is a point $\bar P$ on $\mathcal{D}$ such that 
$P = j(\bar P)$.
Then by Theorem~\ref{Fano}, there exists a skyscraper sheaf $B$ of length 
$1$ on $\mathcal{F}$ supported at $\bar Q = \bar{\psi}(\bar P)$ 
such that $A$ is contained in the subcategory 
generated by the sheaves of the form 
$j_*\bar{\psi}^*B \otimes \mathcal{O}_{\mathcal{X}}
(\sum_{i=1}^{n+1} k_i\mathcal{D}_i)$
for 
\[
\sum_{i=\beta+1}^{n+1} \frac{b_i}{r_i} 
> \sum_{i=1}^{n+1} \frac{a_ik_i}{r_i} \ge 
- \sum_{i=1}^{n+1} \frac{a_i}{r_i}.
\]
Therefore, $A$ contained in $\mathcal{T}$, hence ${}^{\perp}\mathcal{T} = 0$.
\end{proof}

\begin{Cor}\label{cflip}
Assume that $D^b(\text{Coh}(\mathcal{X}^+))$ has a complete exceptional 
collection consisting of sheaves. 
Then so has $D^b(\text{Coh}(\mathcal{X}))$.
\end{Cor}


Department of Mathematical Sciences, University of Tokyo, 

Komaba, Meguro, Tokyo, 153-8914, Japan 

kawamata@ms.u-tokyo.ac.jp

\end{document}